\documentclass[12pt]{amsart}
\usepackage{mathrsfs}
\usepackage{amsfonts}
\usepackage{amssymb}
\usepackage{ifthen}
\usepackage{graphicx}
\usepackage{float}
\usepackage[usenames]{color}

\nonstopmode \numberwithin{equation}{section}
\newtheorem{thm}{Theorem} 
\newtheorem{lem}{Lemma} 
\newtheorem{cor}{Corollary} 

\newtheorem{conj}[equation]{Conjecture}

\theoremstyle{definition}
\newtheorem{defn}{Definition}
\newtheorem{example}{Example}
\newtheorem{prob}[equation]{Problem}
\newtheorem{ques}[equation]{Question}
\newtheorem{exam}[equation]{Example}
\newenvironment{rem}{

\medskip
\noindent \textsl{{\bf Remark. }}}{\medskip}
\newenvironment{rems}{
\bigskip
\noindent \textsl{{\bf Remarks. }}}{\bigskip}

\newcounter {own}
\def\theown {\thesection       .\arabic{own}}

\newenvironment{pf}[1][]{%
 \vskip 3mm
 \noindent
 \ifthenelse{\equal{#1}{}}%
  {{\slshape Proof. }}%
  {{\slshape #1.} }%
 }%
{\qed \smallskip}

\newcounter{alphabet}
\newcounter{tmp}
\newenvironment{Thm}[1][]{\refstepcounter{alphabet}%
\bigskip%
\noindent%
{\bf Theorem \Alph{alphabet}}%
\ifthenelse{\equal{#1}{}}{}{ (#1)}%
{\bf .} \itshape}{\vskip 8pt}

\makeatletter
\newcommand{\Ref}[1]{\@ifundefined{r@#1}{}{\setcounter{tmp}{\ref{#1}}\Alph{tmp}}}
\makeatother

\newcommand{\IN}{{\mathbb N}}
\newcommand{\IC}{{\mathbb C}}
\newcommand{\ID}{{\mathbb D}}

\newcommand{\grad}{{\operatorname{grad}\,}}

\def\be{\begin{equation}}
\def\ee{\end{equation}}

\newcommand{\bee}{\begin{enumerate}}
\newcommand{\eee}{\end{enumerate}}

\newcommand{\blem}{\begin{lem}}
\newcommand{\elem}{\end{lem}}
\newcommand{\bthm}{\begin{thm}}
\newcommand{\ethm}{\end{thm}}
\newcommand{\bcor}{\begin{cor}}
\newcommand{\ecor}{\end{cor}}
\newcommand{\beg}{\begin{exam}}
\newcommand{\eeg}{\end{exam}}
\newcommand{\begs}{\begin{examples}}
\newcommand{\eegs}{\end{examples}}
\newcommand{\bdefe}{\begin{defn}}
\newcommand{\edefe}{\end{defn}}
\newcommand{\bprob}{\begin{prob}}
\newcommand{\eprob}{\end{prob}}
\newcommand{\bques}{\begin{ques}}
\newcommand{\eques}{\end{ques}}
\newcommand{\bei}{\begin{itemize}}
\newcommand{\eei}{\end{itemize}}
\newcommand{\bcon}{\begin{conj}}
\newcommand{\econ}{\end{conj}}
\newcommand{\bcons}{\begin{conjs}}
\newcommand{\econs}{\end{conjs}}
\newcommand{\bprop}{\begin{propo}}
\newcommand{\eprop}{\end{propo}}
\newcommand{\br}{\begin{rem}}
\newcommand{\er}{\end{rem}}
\newcommand{\brs}{\begin{rems}}
\newcommand{\ers}{\end{rems}}
\newcommand{\bo}{\begin{obser}}
\newcommand{\eo}{\end{obser}}
\newcommand{\bos}{\begin{obsers}}
\newcommand{\eos}{\end{obsers}}
\newcommand{\bpf}{\begin{pf}}
\newcommand{\epf}{\end{pf}}
\newcommand{\ba}{\begin{array}}
\newcommand{\ea}{\end{array}}
\newcommand{\beq}{\begin{eqnarray}}
\newcommand{\beqq}{\begin{eqnarray*}}
\newcommand{\eeq}{\end{eqnarray}}
\newcommand{\eeqq}{\end{eqnarray*}}

\newcommand{\ds}{\displaystyle}

\begin{document}
\bibliographystyle{amsplain}
\title {Criteria of univalence and fully {\Large$\alpha$}--accessibility for {\small$p$}--harmonic and
{\small$p$}--analytic functions}

\author{K. F. Amozova, E. G. Ganenkova, and S. Ponnusamy}

\address{K. F. Amozova, Faculty of Mathematics and Information Technology,
Petrozavodsk State University, 33, Lenina st., 185910 Petrozavodsk,
Russia.} \email{amokira@rambler.ru}

\address{E. G. Ganenkova,
Faculty of Mathematics and Information Technology, Petrozavodsk
State University, 33, Lenina st., 185910 Petrozavodsk, Russia.}
\email{g\_ek@inbox.ru}

\address{S. Ponnusamy,
Indian Statistical Institute (ISI), Chennai Centre, SETS (Society
for Electronic Transactions and Security), MGR Knowledge City, CIT
Campus, Taramani, Chennai 600 113, India.
}
\email{samy@isichennai.res.in, samy@iitm.ac.in}

\subjclass[2000]{Primary: 31A05,  31A30,  31C05,   30C45, 52A30; Secondary: 30C20}
\keywords{Univalence, $p$--harmonic , $p$--analytic function,  and fully $\alpha$--accessible function, $\alpha$--accessible domain,
starlike domain and Elliptic PDE}

\begin{abstract}
In this article, we present univalence criteria for polyharmonic and polyanalytic functions.
Our approach yields new a criterion for a polyharmonic functions to be fully $\alpha$--accessible.
Several examples are presented to illustrate the use of these criteria.
\end{abstract}


\maketitle \pagestyle{myheadings}
\markboth{K. F. Amozova, E. G. Ganenkova, and S. Ponnusamy}{Criteria of fully {$\alpha$}--accessibility for {$p$}--harmonic  functions}

\section{Introduction}

Among the most widely studied objects in function theory are the complex-valued $C^1$-solutions $f=u+iv$ of Cauchy-Riemann equation $f_{\overline{z}}(z)=0$
known as analytic functions. The next one is the complex-valued $C^2$-solutions $f$ of the Laplace equation $\Delta f=0$ known as harmonic functions, where
$$\Delta=4\frac{\partial^{2}}{\partial z\partial\overline{z}}=\frac{\partial ^2}{\partial x^2} +
\frac{\partial ^2}{\partial y^2}
$$
stands for the Laplacian (for details see, for example, \cite{Balk,CPW1,Li}).
Nevertheless, there are non-analytic and non-harmonic functions which play significant role in the study of function
spaces and analog properties, for example. Let $\mathbb{D}=\{z\in\mathbb{C}:\,|z|<1\}$
be the unit disk and $f\in C^{2p}(\ID)$, i.e., $f$ is a complex-valued function which is $2p$--times continuously differentiable
on $\ID$. Then $f$ is called $p$--harmonic (or polyharmonic of order $p-1$) in $\ID$ if $f$ satisfies the $p$--harmonic equation $\Delta^p f=0$
in $\ID$, where $p\in \mathbb{N}$.  If $f$ is $p$--harmonic for some $p$, then we simply say that  it is polyharmonic.
Obviously, when $p = 1$ (resp. $p = 2$), $f$ is harmonic (resp. biharmonic) in $\mathbb{D}$.
Also, it is clear that every harmonic mapping is $p$--harmonic for each $p\geq2.$

For $\alpha\in\mathbb{R}$ and $z\in\mathbb{D}$, let
$$T_{\alpha}=-\frac{\alpha^{2}}{4}(1-|z|^{2})^{-\alpha-1}+\frac{\alpha}{2}(1-|z|^{2})^{-\alpha-1}
\left(z\frac{\partial}{\partial z}+\overline{z}\frac{\partial}{\partial
\overline{z}}\right)+\frac{1}{4}(1-|z|^{2})^{-\alpha}\Delta
$$
be the {\it second order elliptic partial differential operator}. Of particular interest to our analysis is the homogeneous partial
differential equation
\be\label{eq-CV}
T_{\alpha}(f)=0~\mbox{in}~\mathbb{D}
\ee
and its associated {\it Dirichlet boundary value problem}:
\be\label{eq-1}
\begin{cases}
\displaystyle T_{\alpha}(f)=0
& \mbox{in } \mathbb{D},\\
\displaystyle~ f=f^{\ast} &\mbox{on } \partial\mathbb{D}.
\end{cases}
\ee
Here, the boundary data $f^{\ast}\in\mathfrak{D}'(\partial\mathbb{D})$ is a
{\it distribution} on the boundary  $\partial\mathbb{D}$ of $\mathbb{D}$,
and the boundary condition in (\ref{eq-1}) is interpreted in the
distributional sense that $f_{r}\rightarrow f^{\ast}$ in
$\mathfrak{D}'(\partial\mathbb{D})$ as $r\rightarrow1-$, where
\be\label{eq-0.1}
f_{r}(e^{i\theta})=f(re^{i\theta}),~e^{i\theta}\in\partial\mathbb{D},
\ee
for $r\in[0,1)$. In \cite{Olo-2014}, Olofsson proved  that, for parameter values
$\alpha>-1$, if a function $f\in C^{2}(\mathbb{D})$ satisfies (\ref{eq-CV}) with
$\lim_{r\rightarrow1-}f_{r}=f^{\ast}\in\mathfrak{D}'(\partial\mathbb{D})$,
then  it has the form of a {\it Poisson type integral}
\be\label{eq-0.2}
f(z)=\frac{1}{2\pi}\int_{0}^{2\pi}K_{\alpha}(ze^{-i\tau})f^{\ast}(e^{i\tau})\,d\tau~\mbox{  $z\in\mathbb{D}$,}
\ee
where
$$K_{\alpha}(z)=c_{\alpha}\frac{(1-|z|^{2})^{\alpha+1}}{|1-z|^{\alpha+2}},
$$
$c_{\alpha}=\big(\Gamma(\alpha/2+1)\big)^{2}/\Gamma(1+\alpha)$ and $\Gamma(s)=\int_{0}^{\infty}t^{s-1}e^{-t}\,dt$
for $s>0$ is the  standard Gamma function (see also \cite{AH-2014,CV-2015}).

If we take $\alpha =2(p-1)$, then $f\in C^{2p}(\ID)$ is $p$-harmonic. (cf. \cite{AH-2014,CPW1,CPW2-2013,CPW3-2012,CV-2015}).
Furthermore, Borichev and Hedenmalm \cite{AH-2014} proved that
$$(1-|z|^{2})^{p}\Delta^{p}=4(1-|z|^{2})T_{0}\circ4(1-|z|^{2})^{2}T_{2}\circ\cdots\circ4(1-|z|^{2})^{p}T_{2(p-1)}.
$$
In particular, if $\alpha=0$, then $f$ is harmonic (cf. \cite{Du}).

One of the aims of this article is to study polyharmonic and polyanalytic functions defined on $\ID$.
There are a lot of articles devoted to univalence conditions for
analytic functions in  $\mathbb{D}$, but there exists  only a few in the case of
harmonic functions. See \cite{Clunie-Small-84,Li,PQ-2013} and the references therein.
One of them was obtained recently by Starkov in
\cite{Star} in terms of the coefficients of analytic and co-analytic
parts of the harmonic function. This criterion was established by using
the method of Bazilevich who obtained analogous theorem for analytic functions.

\begin{Thm} (cf. \cite{Baz})\label{ThA}
An analytic function $F(z)=\sum_{n=1}^\infty c_nz^n$ in $\mathbb{D}$ is univalent in $\mathbb{D}$ if and only if  for each
$z\in\mathbb{D}$ and each $t\in[0,\pi/2]$,
\begin{equation}\label{eq0000}\sum_{n=1}^\infty c_n \frac {\sin nt} {\sin t}z^{n-1}\ne 0,
\quad \left(\frac {\sin nt} {\sin t}\right)\Big\vert_{t=0}=n.
\end{equation}
\end{Thm}

In this article, we first present univalence criterion for $p$--harmonic functions in $\mathbb{D}$ (Theorem \ref{th0}) and
demonstrate its use by presenting a number of interesting and simple examples. As a consequence, we also state a univalence criterion for $p$--analytic
functions (Theorem \ref{th0a}). In Section \ref{sec2}, we include the definition and some related developments concerning $\alpha$-accessible domains, and
then present necessary and sufficient conditions for $p$--harmonic functions to be
$\alpha$-accessible, which, in particular, produces a criterion for fully starlike functions.

\subsection{Univalence criterion for polyharmonic functions}

%
%

We begin to recall that $F\in C^{2p}(\ID)$ is $p$--harmonic in $\mathbb{D}$ if and only if $F$ has the
following representation (see \cite{CPW1}):
$$F(z)=\sum\limits_{k=1}^p|z|^{2(k-1)}F_{p-k+1}(z),
$$
where each $F_{p-k+1}=h_{p-k+1}+\overline{g}_{p-k+1},$ $k=1,\dots,p,$ is harmonic in $\ID$. Without loss of generality we may assume that
$h_{p-k+1}(0)=0=g_{p-k+1}(0)$ for convenience. Recall that every harmonic function $f$ in $\ID$ can be written as
$f=h+\overline{g}$ where $h$ and $g$ are analytic with $g(0)=0$. Then
$J_f= |f_z|^2-|f_{\overline{z}}|^2=|h'|^2-|g'|^2$ denotes the Jacobian of $f$. We say that  $f$ is a sense-preserving harmonic mapping
if $J_f(z)>0$ in $\ID$.

\bthm \label{th0}
Let $F(z)=\sum_{k=1}^p|z|^{2(k-1)}F_{p-k+1}(z)$ be $p$-harmonic in $\mathbb{D}$ and univalent in a neighborhood of the origin, where
$$F_{p-k+1}(z)=\sum_{n=1}^{\infty} \left (a_n^{(p-k+1)} z\thinspace^n+b_n^{(p-k+1)}{\overline{z}}\thinspace^n\right ) ~\mbox{ ($k=1,\dots,p$)}
$$
are harmonic in $\mathbb{D}.$ Then the function $F(z)$ is univalent in $\mathbb{D}$ if and only if for each
$z\in\mathbb{D}\backslash\{0\}$ and $t\in\left(0,\pi/2\right]$ the following condition holds:
\begin{equation}\label{eq0}
\sum\limits_{k=1}^p|z|^{2(k-1)}\sum\limits_{n=1}^{\infty}\left(a_n^{(p-k+1)} z\thinspace^n-b_n^{(p-k+1)}{\overline{z}}\thinspace^n\right)
\frac{\sin nt}{\sin t}\neq0.
\end{equation}
\ethm
\bpf
We use the method of proof of Bazilevich.

\vspace{6pt}
\noindent{\underbar{Necessity}.} Let $F(z)$ be univalent in
$\mathbb{D}$ and take $z_1=re^{i\theta_1}, $
$z_2=re^{i\theta_2}\in\mathbb{D},$ $z_1\neq z_2,$
$\theta_1<\theta_2\leq\theta_1+\pi,$ $r\in(0, 1).$ Then we have
\be\label{eq-new1}
\frac{F(z_2)-F(z_1)}{z_2-z_1}\neq 0
\ee
which may be rewritten as
\begin{equation}\label{1}
\sum_{k=1}^p r^{2(k-1)}\sum_{n=1}^\infty \left(a_n^{(p-k+1)}r^{n-1}\frac {e^{in\theta_2}-e^{in\theta_1}} {e^{i\theta_2}-e^{i\theta_1}} + b_{n}^{(p-k+1)}r^{n-1}\frac {e^{-in\theta_2}-e^{-in\theta_1}} {e^{i\theta_2}-e^{i\theta_1}}\right)\ne 0.
\end{equation}
Since
$e^{in\theta_2}-e^{in\theta_1}=e^{in(\theta_1+\theta_2)/2}\left ( e^{in(\theta_2-\theta_1)/2}-e^{-in(\theta_2-\theta_1)/2}\right )$
for $n\in \{\pm 1, \pm2, \ldots \}$, it follows easily that
$$\frac {e^{in\theta_2}-e^{in\theta_1}} {e^{i\theta_2}-e^{i\theta_1}} =e^{iT(n-1)}\frac{\sin nt}{\sin t},\quad
t=\frac{\theta_2-\theta_1}{2}\in(0,\pi/2], ~~ T=\frac{\theta_1+\theta_2}{2},
$$
and thus, replacing $n$ by $-n$ gives
$$\frac{e^{-i n\theta_2}-e^{-i n\theta_1}}{e^{i\theta_2}-e^{i\theta_1}}= -e^{-i T (n+1)}\frac{\sin nt}{\sin t}  .
$$
Therefore, (\ref{1}) takes the form
\begin{equation}\label{eq001}
\sum\limits_{k=1}^pr^{2(k-1)}\sum\limits_{n=1}^{\infty}r^{n-1}\left(a_n^{(p-k+1)} e^{i T (n-1)}-
b_n^{(p-k+1)}e^{-i T (n+1)}\right)\frac{\sin nt}{\sin t}\neq0,
\end{equation}
or equivalently,
\begin{equation}\label{eq00}
e^{-i T}\left (\sum\limits_{k=1}^p|z|^{2(k-1)}\sum\limits_{n=1}^{\infty}\left(a_n^{(p-k+1)} z^{n}-
b_n^{(p-k+1)}{\overline{z}}^{n}  \right)
\frac{\sin nt}{\sin t}\right )\neq0,
\end{equation}
where $z=re^{i T}\neq0,$ $t\in\left(0, \pi/2\right].$ The last relation is clearly equivalent to the condition \eqref{eq0} and the proof of the
necessary part is complete.

\vspace{6pt}
\noindent{\underbar{Sufficiency}.} Let $F(z)$ satisfy the condition (\ref{eq0}). As shown above, this is equivalent to the condition
$$ \frac{F(z_2)-F(z_1)}{z_2-z_1}\neq0
$$
for each $z_1=re^{i\theta_1}\neq z_2=re^{i\theta_2},$ $r\in(0, 1),$ $\theta_1,\theta_2\in\mathbb{R.}$
This means that $F(z)$ is univalent on every circle $\{z\in\mathbb{C}:\, |z|=r\}.$

Suppose on the contrary that $F(z)$ is not univalent in $\mathbb{D}.$ Since  $F(z)$ is univalent in a neighborhood of the origin, by assumption,
there exists a disk $\rho\mathbb{D}:=\mathbb{D}_\rho=\{z\in\mathbb{C}:\, |z|<\rho \}$, $\rho<1$, such that
$F(z)$ is univalent in $\mathbb{D}_\rho$ and we can find points $z_1,
z_2\in\partial\mathbb{D}_\rho=\{z\in{\mathbb{C}}:\, |z|=\rho\},$ $z_1\neq z_2,$ for which
$F(z_1)=F(z_2).$ So $F(z)$ is not univalent on the circle $\partial\mathbb{D}_\rho.$
The contradiction proves the sufficiency.
\epf

\br
Univalence condition (\ref{eq001}) for $p$--harmonic functions is an analog condition (\ref{eq0000}) for analytic functions with
$z=re^{iT}.$ But, unlike the analytical case, the non-vanishing condition (\ref{eq001}) is not necessarily true for $z=0.$
If $F(z)$ is analytic and we take $z=0$ in (\ref{eq0000}), then
we obtain  $c_1=F'(0)=J_F(0)\ne0$  and thus, $F$ is univalent in a neighborhood of the origin. But for $p$--harmonic functions,
it may happen that $F$ is univalent even if $J_F(0)=0.$ For example, for the biharmonic function $F(z)=|z|^2z$, we have $J_F(z)=3|z|^4.$
For such functions the left-hand side of the condition (\ref{eq001}) at  $z=0$ equals $J_F(0).$ Therefore, (\ref{eq001}) does not hold for $z=0.$
\er

\subsection{Examples of polyharmonic univalent functions}
We now demonstrate by examples how Theorem \ref{th0} can be applied in practice. To do this, we consider for instance
$p$--harmonic functions of the form
\be\label{eg2-ext1}
F(z)=|{z}|^{2(p-1)}G(z)+K(z),
\ee
where $p\geq 2$, $G$ and $K$ are harmonic in $\ID$.  Moreover, for $F$ of the form \eqref{eg2-ext1}, we have
\beq\label{eg2-ext2}
J_F(z)& =& |F_z(z)|^2-|F_{\overline z}(z)|^2\nonumber\\
& =& \big |(p-1)|z|^{2(p-2)}\overline z G(z)+|z|^{2(p-1)} G_{z}(z)+K_z(z)\big |^2 \nonumber\\
&& ~~~~ -\big|(p-1)|z|^{2(p-2)}z G(z)+|z|^{2(p-1)} G_{\overline z}(z)+K_{\overline z}(z)\big|^2
\eeq
and observe that $J_F(0)=J_K(0)$. For our purpose, in the examples below, we consider some special choices for $G$ and $K$.

\begin{example}\label{eg2}
In \eqref{eg2-ext1}, we let  $G(z)=z$ and $K(z)=z+\lambda {\overline z}^n$,  $|\lambda|\leq 1/n$ with $n\geq 2$.
Clearly, $K$ is harmonic and univalent in $\ID$ such that $J_K(0)=1$. Moreover, by \eqref{eg2-ext2}, we obtain that $J_F(0)=J_K(0)=1$
and thus, $F(z)$ is univalent in a neighborhood of the origin.

Next to prove that $F$ is univalent in $\ID$, it suffices to show that the
condition \eqref{eq0} holds. In our case this condition is equivalent to
\be\label{eq-new3}
z\left(|z|^{2(p-1)} +1\right )\neq \lambda \overline{z}\thinspace ^n\cdot\frac{\sin nt}{\sin t}
\ee
for all $z\in\mathbb D\backslash\{0\}$ and $t\in\left(0, \pi/2\right]$. We show this by a contradiction.
Suppose that \eqref{eq-new3} is not true. Then
$$z\left(|z|^{2(p-1)}+1\right)=\overline{z}\left(\lambda \overline
z\thinspace ^{n-1} \cdot\frac{\sin nt}{\sin t}\right)
$$
for some $z\in\mathbb D\backslash\{0\}$ and for some $t\in\left(0, \pi/2\right].$ This yields
\be\label{eq-new4}
|z|^{2(p-1)}+1=\left|\lambda  \overline{z}\thinspace
^{n-1}\cdot\frac{\sin nt}{\sin t}\right|
\ee
at such a point $z$ and some $t$ in that interval. Clearly, $ |z|^{2(p-1)}+1>1$ in $\mathbb D\backslash\{0\}.$
%
On the other hand, since $|\sin nt| \le n |\sin t|$  for all $t \in [0, \pi/2]$ and $n \in \IN$ (which may easily be verified by a method of induction),
it follows that for any $z\in\mathbb D\backslash\{0\}$ and for any $t\in\left(0, \pi/2\right]$,
$$\left|\lambda  \overline{z}\thinspace ^{n-1}\cdot\frac{\sin nt}{\sin
t}\right|<|\lambda| \left|\frac{\sin nt}{\sin t}\right|\leq |\lambda| n\leq 1
$$
which contradicts \eqref{eq-new4}. Therefore, $F(z)$ is univalent in $\mathbb D.$ The graph of $F(z)$ with $\lambda =1/n$ and
for certain values of $p$ and $n$ are shown in Figures \ref{fig2ab} and \ref{fig2cd}.

\begin{figure}[H]
\begin{center}
\includegraphics[height=5.0cm, width=5.5cm, scale=1.3]{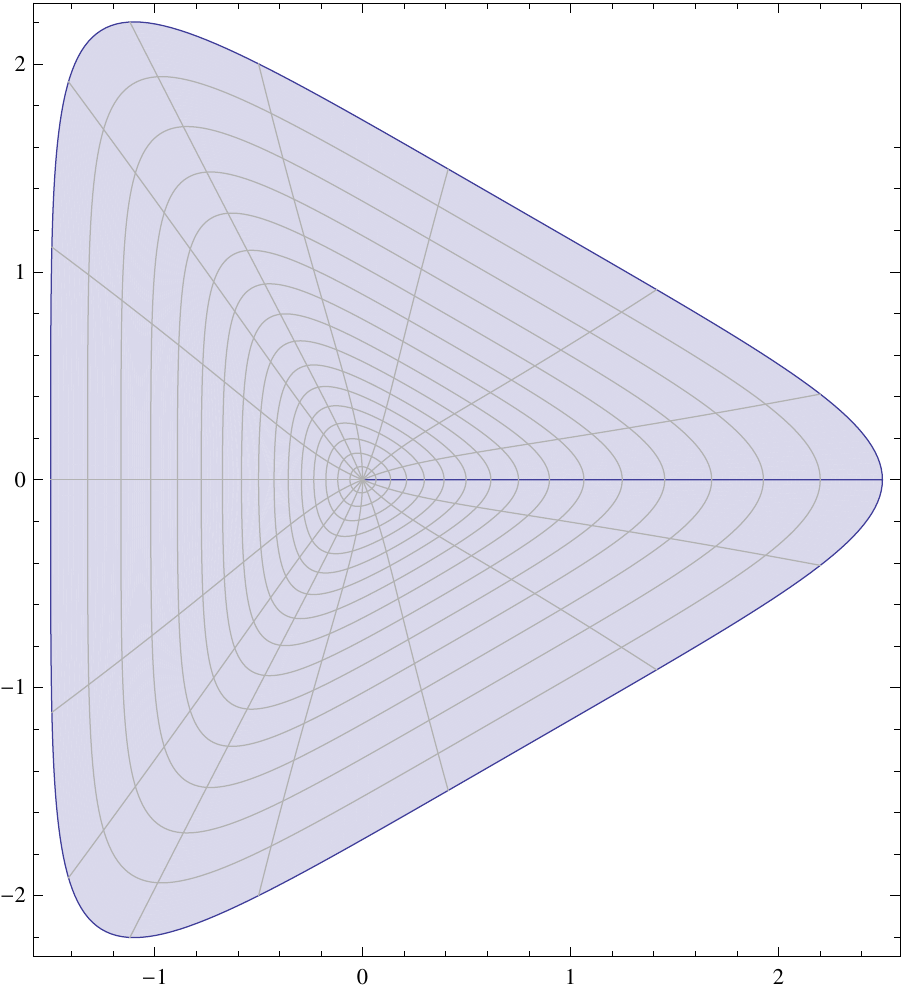}
\hspace{1cm}
\includegraphics[height=5.0cm, width=5.5cm, scale=1.3]{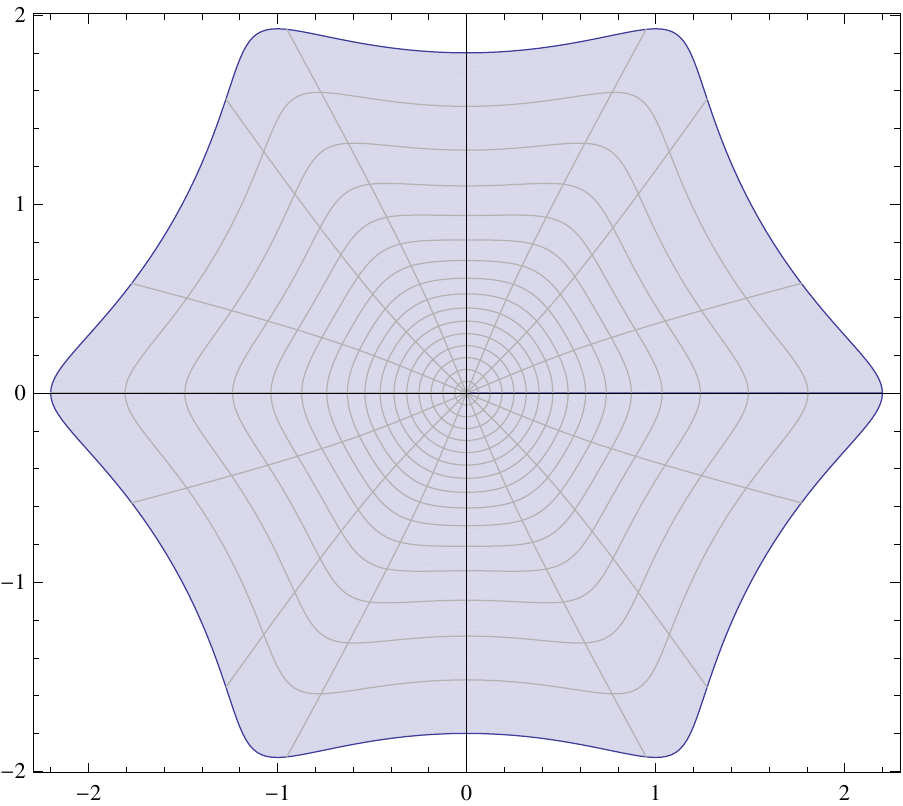}
\end{center}
(1) $n=2$, $p=2$ \hspace{4.3cm} (2) $n=5$, $p=3$
\caption{The range of $F(z)=|z|^{2(p-1)}z+ z+(1/n) (\overline{z})^n$ for certain values of $p$ and $n$
} \label{fig2ab}
\end{figure}


\begin{figure}[H]
\begin{center}
\includegraphics[height=5.0cm, width=5.5cm, scale=1.3]{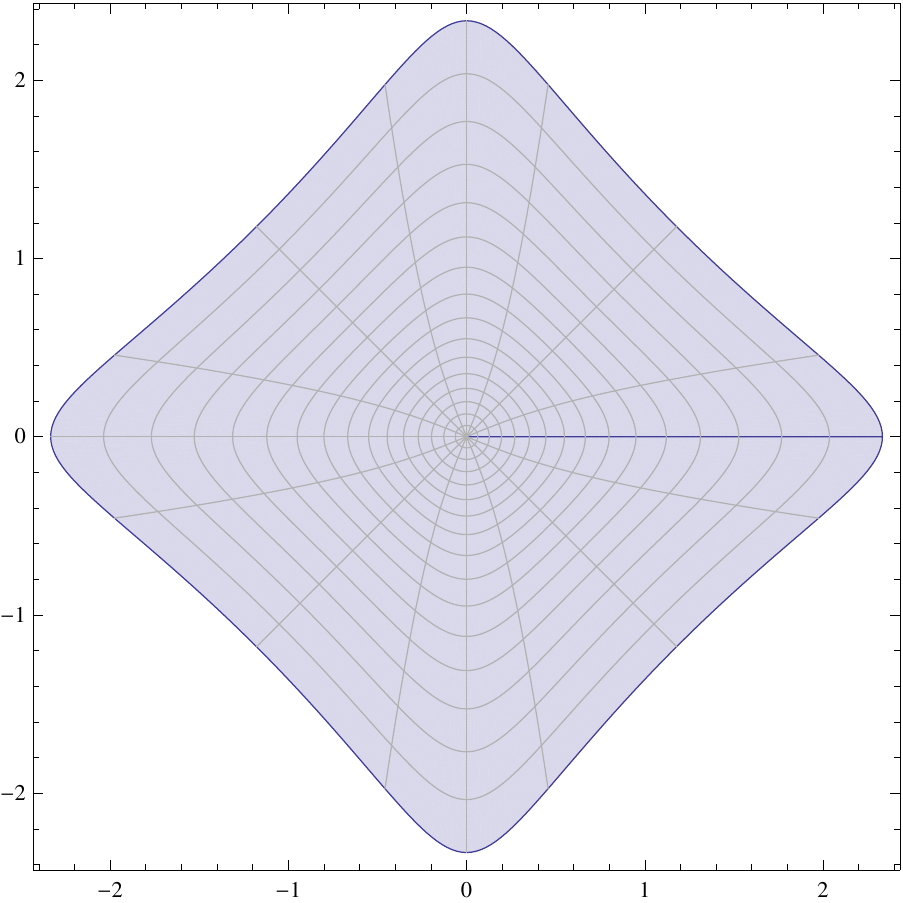}
\hspace{1cm}
\includegraphics[height=5.0cm, width=5.5cm, scale=1.3]{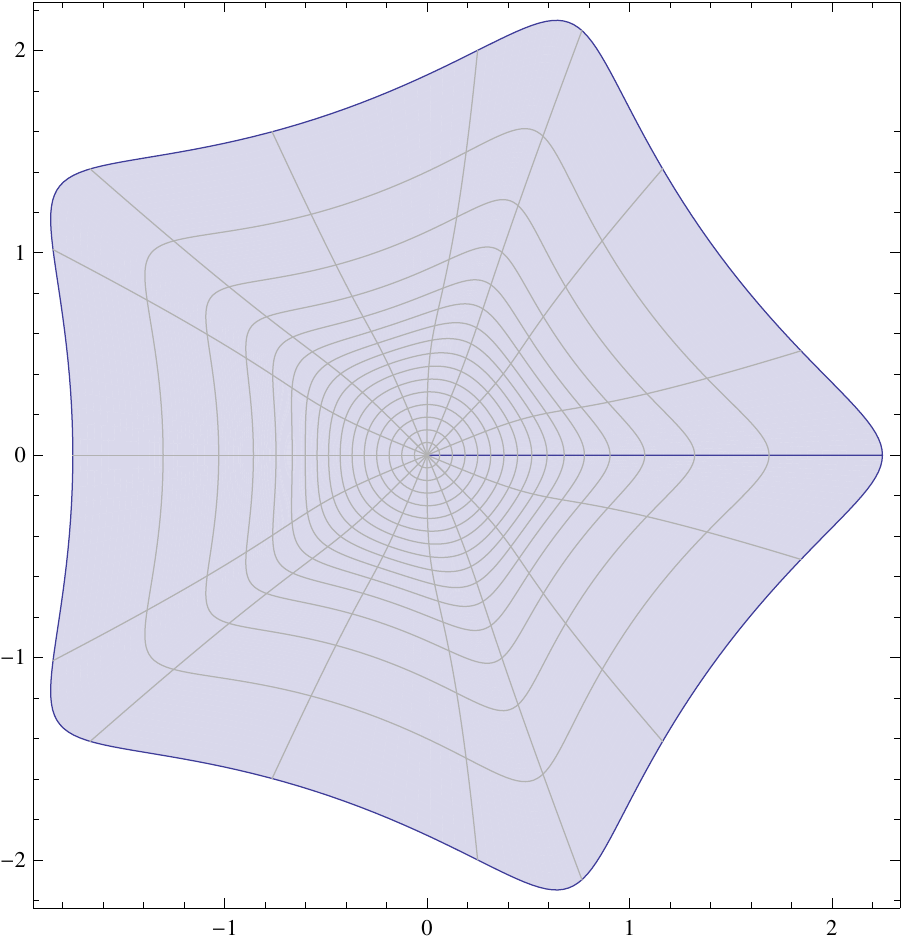}
\end{center}
(1) $n=3$, $p=2$ \hspace{4.3cm} (2) $n=4$, $p=5$

\caption{The range of $F(z)=|z|^{2(p-1)}z+ z+(1/n) (\overline{z})^n$ for certain values of $n$ and $p$
}\label{fig2cd}
\end{figure}
\end{example}

\begin{example}\label{eg1}
In \eqref{eg2-ext1}, we let $G(z)=z$ and
$$K(z)=a\left(\displaystyle\frac{z+c}{1+\overline{c}z}\right)+b\overline{z},
$$
where  $a,b$ and $c$ are complex numbers such that $c\in\mathbb{D},$  $|b|\ne|a| (1-|c|^2)$ and
$$|a|\geq (1+|b|)\left (\frac{1+|c|}{1-|c|}\right ).
$$
Then, by \eqref{eg2-ext1} and \eqref{eg2-ext2}, $J_F(0)=J_K(0)=\left|a(1-|c|^2)\right|^2-|b|^2\ne0$. Therefore,
the function $F(z)=|z|^{2(p-1)}G(z)+K(z)$ is univalent in a neighborhood of the origin.
Next, we show that condition \eqref{eq0} holds for our function $F(z).$ Since
$$\frac{z+c}{1+\overline{c}z}=c+(1-|c|^2)\sum\limits_{n=1}^{\infty}(-\overline{c})^{n-1}z^n,
$$
condition  \eqref{eq0} takes the form
$$
|z|^{2(p-1)}z+ \frac{a(1-|c|^2)}{\sin t}\sum\limits_{n=1}^{\infty}\left ( (-\overline{c})^{n-1} \thinspace\sin nt \right )z^n \ne b\overline{z}
$$
for $z\in\mathbb D\backslash\{0\}$ and $t\in\left(0, \pi/2\right].$

%
Using the identity $\sin nt=(e^{int}-e^{-int})/(2i)$, we can simplify the last expression in an equivalent form as
\begin{equation}\label{equat000}
b\overline{z} \ne z\left ( |z|^{2(p-1)}+ \displaystyle\frac{a(1-|c|^2)}{1+2\overline{c} z \cos t +{\overline{c}}^2z^2}\right )
\end{equation}
and thus, the proof is complete if we show that \eqref{equat000} holds for all $z\in\mathbb D\backslash\{0\}$ and $t\in\left(0, \pi/2\right].$

We complete the proof by a contradiction. Suppose that \eqref{equat000} is not true. Then there exists an $z\in\mathbb D\backslash\{0\}$ and a $t\in\left(0, \pi/2\right]$
such that
$$  b\overline{z}=z\left(|z|^{2(p-1)}+\displaystyle\frac{a(1-|c|^2)}{1+2\overline{c} z \cos t +{\overline{c}}^2z^2}\right).
$$
But then,
\beqq
|b|&=&\left ||z|^{2(p-1)}+\displaystyle\frac{a(1-|c|^2)}{1+2\overline{c} z \cos t +{\overline{c}}^2z^2}\right |\\
&\geq& \frac{|a|\, (1-|c|^2)}{|1+2\overline{c}z \cos t +{\overline{c}}^2z^2|}-|z|^{2(p-1)}\\
&>&|a|\left (\displaystyle\frac{1-|c|}{1+|c|}\right )-1,
\eeqq
which contradicts our initial assumptions on the parameters $a,b$ and $c.$ Consequently, $F$ is univalent in $\mathbb{D}.$
\end{example}


Let us now present an unbounded $p$--harmonic univalent function in $\mathbb{D}$.

\begin{example}\label{eg3}
For $p\geq2$ and $\mu\in\IC$ such that $0<|\mu|\le 1/2$, we consider
\be\label{eg3-eq1}
F(z)=|z|^{2(p-1)}\mu\overline{z}+\frac{1+z}{1-z}= |z|^{2(p-1)}\mu\overline{z}+ 1+2\sum\limits_{n=1}^{\infty}z^n.
\ee
We see that
$$J_F(z) =\left |\mu (p-1)|z|^{2(p-2)}\overline{z}^2 +\frac{2}{(1-z)^2}\right |^2 -\left |\mu (p-1)|z|^{2(p-1)} \right |^2
$$
and in particular, $J_F(0)=J_K(0)=2.$ Therefore, $F$ is univalent in a neighborhood of the origin. In this case, as in the previous example,
 condition \eqref{eq0} takes the form
$$ 0\neq -|z|^{2(p-1)}\mu\overline{z}+2\sum\limits_{n=1}^{\infty} \frac{\sin nt}{\sin t} z^n = -|z|^{2(p-1)}\mu\overline{z}+\frac{2z}{1-2z\cos t +z^2}
$$
for $z\in\mathbb D\backslash\{0\}$ and $t\in\left(0, \pi/2\right].$ Finally, since
$$ |z|^{2(p-1)}|\mu\overline{z}|< \frac{|z|}{2} ~\mbox{ and }~ \frac{2|z|}{|1-2z\cos t +z^2|}> \frac{|z|}{2} ,
$$
it follows that the last condition obviously holds  for all $z\in\mathbb D\backslash\{0\}$ and $t\in\left(0, \pi/2\right].$
According to Theorem \ref{th0}, the function $F$ defined by \eqref{eg3-eq1} is univalent in $\ID$. The image of $\ID$ under $F$
for (1) $p=2$ and $\mu=1/4$; and (2) $p=5$ and $\mu=1/2$, are shown in Figure \ref{F-unbound}.


\begin{figure}[H]
\begin{center}
\includegraphics[height=5.5cm, width=5.5cm, scale=1.3]{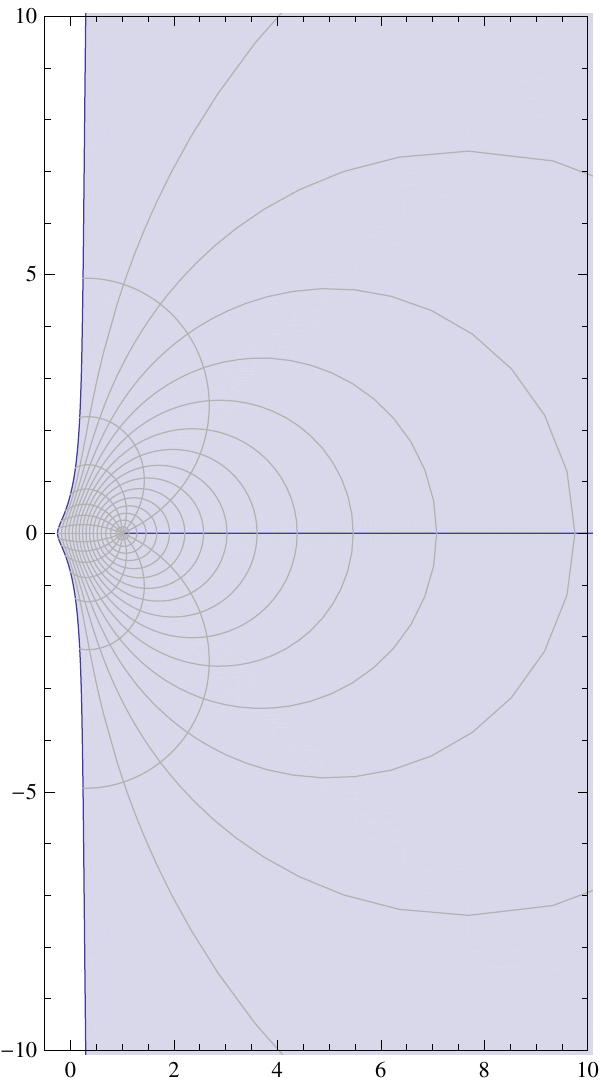}
\hspace{1cm}
\includegraphics[height=5.5cm, width=5.5cm, scale=1.3]{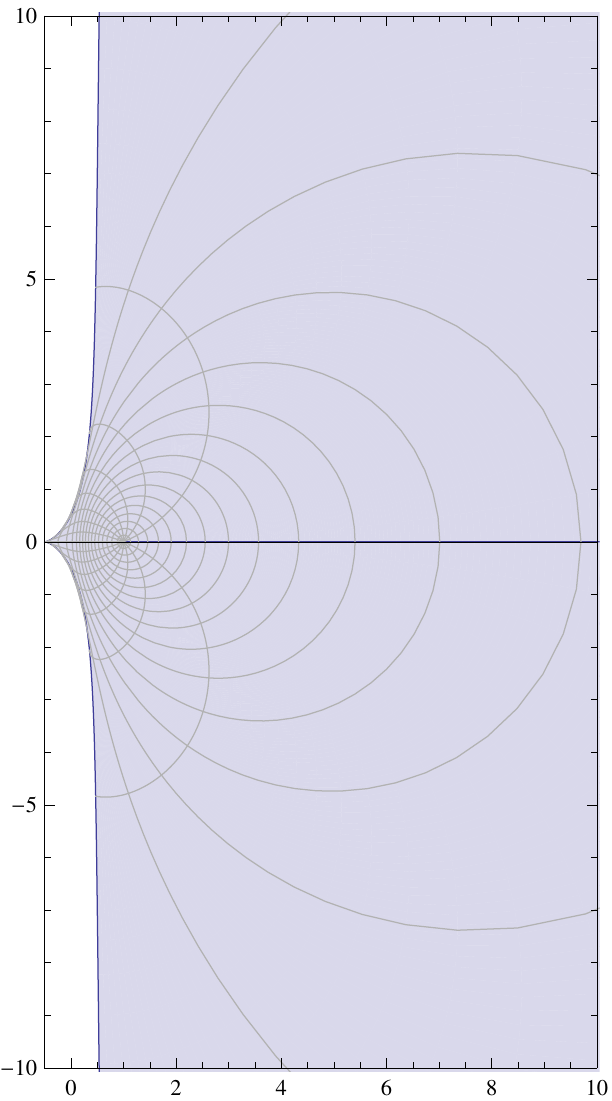}
\end{center}
(1) $p=2$, $\mu=1/4$ \hspace{4.3cm} (1) $p=5$, $\mu=1/2$

\caption{The range of $F(z)=|z|^{2(p-1)}\mu\overline{z} + \frac{1+z}{1-z}$ for certain values of $p$ and $\mu$\label{F-unbound}}
\end{figure}
\end{example}

\subsection{Comparison with a known univalence condition of biharmonic mappings}
At this place it is worth recalling the following sufficient univalence condition for biharmonic mappings (see \cite{Abd}) and to compare it with
Theorem \ref{th0}.

\begin{Thm}\label{th-Abd1}
Let $F(z)=|z|^{2}G(z)+K(z)$ be a biharmonic $(p=2)$
function in the unit disk $\mathbb{D}.$ If $K$ is univalent and $K(\mathbb{D})$ is
a Lavrentiev domain with constant $M,$ i.e. there exists a constant
$0<M<\infty$ such that any two points $w_1, w_2\in K(\mathbb{D})$
can be joined by a path $\gamma\subset K(\mathbb{D})$ of length
$l(\gamma)\leq M|w_1-w_2|,$ and if
\begin{equation}\label{eq}
\displaystyle\frac{2|G|+|G_z|+|G_{\overline z}|}{|K_z|-|K_{\overline z}|}<\displaystyle\frac{1}{M} 
\end{equation}
for $z\in\ID$, then $F$ is univalent in $\mathbb{D}.$
\end{Thm}

We now show that Theorem \Ref{th-Abd1} is not applicable to our earlier examples. For instance in Example \ref{eg2},
we have $G(z)=z$ and  $K(z)=z+\displaystyle\frac{{\overline z}\thinspace^n}{n}$ (with $\lambda =1/n$ and $n\geq 2$).
Clearly, $K(\mathbb{D})$ is a Lavrentiev domain with some constant $M.$ Then the condition \eqref{eq} is equivalent to
$$\displaystyle\frac{2|z|+1}{1-|z|^{n-1}}<\frac1M,
$$
which is clearly not valid when $|z|$ approaches $1^{-}$, because the left-hand
side is unbounded.

In Example \ref{eg1}, put $c=0.$ Then we have
$ G(z)=z$ and $K(z)=az+b\overline{z},$ where $|a|\ge |b| +1.$ The function $K$ is univalent in $\mathbb{D}$ and $K(\mathbb{D})$ is a convex domain. Therefore,
$K(\mathbb{D})$ is a Lavrentiev domain with $M=1.$ Moreover, the condition \eqref{eq} reduces to the form
$$2|z|+1<|a|-|b|,
$$
which is obviously not true for all points in $\ID$ and for any pair $(a,b)$ of complex numbers $a$ and $b$ such that $1\le |a|-|b|<3$.
Thus,  Theorem \Ref{th-Abd1} is again not applicable in this example too. However,  according to Example \ref{eg1} with $c=0$, we obtain that the function
$F(z)=|z|^{2(p-1)}z+az +b\overline{z}$ is univalent in $\ID$ if $p\geq 2$ and $a,b\in\IC$ such that
$|a|\geq |b|+1$.

Finally, in Example \ref{eg3}, we take $\mu=1/2.$ Then $G(z)=\overline{z}/2,$ and $K(z)=(1+z)/(1-z)$ which is a univalent mapping of $\mathbb{D}$
onto the right half-plane ${\rm Re}\, w>0$. Therefore, the condition \eqref{eq} takes the following form
$$|1-z|^2(2|z|+1)<4.
$$
For $z\to-1$, the left-hand side of the inequality tends to $12$ and thus, the last inequality cannot be not true for all $z\in\ID$. Thus, the univalency of
$$F(z)=(1/2)|z|^{2} \overline{z}+ \frac{1+z}{1-z}
$$
follows from Theorem \ref{th0} but not from Theorem \Ref{th-Abd1}.

\subsection{Univalence criterion for polyanalytic functions}

Using the method of proof of Theorem \ref{th0}, we can obtain a univalence criterion for $p$--analytic functions.

Polyanalytic functions $F$ of order $p$ (or simply polyanalytic functions) in a simply connected domain  $D$ are polynomials
in the variable $\overline{z}$ with analytic functions as its coefficients:
$$F(z)=\sum\limits_{k=0}^{p-1}\overline{z}\thinspace^k a_{k}(z),
$$
where $\{a_{k}(z)\}_{k=1}^{p-1}$ are analytic in $D$. Clearly, they
may be also defined as $C^p(D)$-solutions of the generalized
Cauchy--Riemann condition $\partial ^p F/\partial \overline{z}^p =
0$ in $D$ (the Cauchy-Riemann equation of order $p$) (for details
see, for example, \cite{Balk}).
It is not difficult to see that every $p$--analytic function is
$p$--harmonic. But there exist $p$--harmonic functions that are not
$p$--analytic. For example, consider the $p$--harmonic function
$F(z)=|z|^{2(p-1)}(z+\overline z).$ Since $\partial ^p F/\partial
\overline{z}^p =p!z^{p-1}\not \equiv 0$ in $\ID,$ it follows that $F(z)$ is not
$p$--analytic.

Our interest is when $D$ is the unit disk $\ID$. Obviously, some basic properties of analytic functions ceases to be true for polyanalytic function.
For example, $F(z)=1-|z|^2$ is polyanalytic of order $2$ but is not analytic. Moreover, polyanalytic functions are closely related to the usual
analytic functions ($p=1$ case). A polyanalytic function of order $p=2$ is
called bianalytic and is of the form $a_0(z)+\overline{z}a_1(z)$. The functions of this form
can be an efficient tool for solving problems in the planar theory of elasticity. In view of this reasoning, it might be useful to state the
analog of Theorem \ref{th0} for polyanalytic functions. Since the proof is similar we omit the details.

\bthm\label{th0a}
Let $F(z)=\sum_{k=0}^{p-1}\overline{z}\thinspace^k a_{k}(z)$
be a $p$--analytic in $\mathbb{D}$ and univalent in a neighborhood
of the origin, where
$$a_{k}(z)=\sum\limits_{n=0}^{\infty} c_n^{(k)}z^n, \,\,\, k=1,\dots,p-1,
$$
are analytic functions in $\mathbb{D}.$ Then the function $F(z)$ is univalent in $\mathbb{D}$ if and only if for each
$z\in\mathbb{D}\backslash\{0\}$ and $t\in\left(0, \pi/2\right]$ the following condition holds
$$ \sum\limits_{k=0}^{p-1}\overline{z}^{k}\sum\limits_{n=0}^{\infty}c_n^{(k)}  \frac{\sin (n-k)t}{\sin t}\, z^n \neq0.
$$
\ethm

\section{Fully $\alpha$-accessible p--harmonic functions}\label{sec2}
In some fields of mathematics such as the theo\-ry of integral
representations of functions, imbedding theorems, the questions of
the boundary behavior of functions, and the sol\-va\-bi\-li\-ty of
Dirichlet problem (cf. \cite{Adam, Bes, Dol, Zar}) it is important for a domain $\Omega\subset{\Bbb
R}^n$ of a function to possess the cone condition. This means that
there exist numbers $\alpha\in[0,1]$ and $H\in(0,\infty)$ such that
for all points $p\in\Omega$ the right circle cone $V(l(p),H)$ with
vertex $p$ of fixed opening $\alpha\pi,$ height $H,$ and axis vector
$l(p)$ is contained in $\Omega$ (see \cite{math}).

In \cite{LSt} the notion of $\alpha$--accessible domain was defined.
It was shown that such domains posses the cone condition with
$l(p)=-p.$

\bdefe (cf. \cite{LSt}) Let $\alpha\in[0, 1].$ A domain
$\Omega\subset{\Bbb R}^n$ is said to be $\alpha$--accessible with
respect to a point $a\in \Omega$ if, for each finite point $p$ on the boundary $\partial \Omega$,
there exists a number $r=r(p)>0$ such that the cone
{\small$$K_{+}(p,a,\alpha,r)=\left\{x \in {\mathbb{R}\thinspace}^n[p,r]:
\left \langle x-p,\displaystyle\frac{p-a}{\|p-a\|}\right\rangle \geq
\|x-p\|\cos\left (\frac{\alpha\pi}{2}\right ), \, \|x-p\|\leq
r\right\}$$}
is contained in ${\Bbb R}^n\backslash \Omega,$ where $\langle . , .\rangle$ denotes the usual inner product in ${\Bbb R}^n$.
\edefe

Detailed information about $\alpha$--accessible domains can be also found in \cite{Am1, Am2, AmGan, AmGanPiter, AmSt, An}.
In \cite{AmGanPiter}, it was shown that $1$-accessible domains are only disks centered at $a.$
In \cite{LSt}, it was observed that $\alpha$--accessible domains are
also starlike. In the case of smooth boundary, the following criterion of $\alpha$--accessibility  was obtained.

\begin{Thm} {\rm (cf. \cite{LSt})}\label{ThB} Let a domain $\Omega\subset{\Bbb
R}^n,$ $0\in\Omega,$ and $\Omega$ be determined by the following conditions:
\begin{enumerate}
\item $\partial\Omega$ is a smooth real manifold of
dimension $n-1,$ determined by the equation $F(x) = 0,$ $x\in{\Bbb
R}^n,$
where $F$ is a smooth function;
\item domain $\Omega$ is given by the inequality $F(x) <
0.$
\end{enumerate}
Then $\Omega$ is an $\alpha$--accessible domain,
$\alpha\in[0,1),$ if and only if for every point $p\in\partial\Omega$,
$$\left  \langle \frac{\grad F(p)}{\|\grad F(p)\|}, \frac{p}{\|p\|}\right  \rangle \geq
\sin\displaystyle\frac{\alpha\pi}{2}.
$$
\end{Thm}

In the planar case ($n=2$) such domains were studied by
Stankiewicz \cite{St1, St2} and also by  Brannan and Kirwan \cite{BrK}. In
\cite{St1, St2},  Stankiewicz considered the class ${\mathcal
S}_{1-\alpha}$ of all $(1-\alpha)$--starlike analytic functions
$f(z)=a_1z+a_2z^2+\dots,$ $a_1\neq 0,$ in $\mathbb D$ such that
$$\Big|\arg\dfrac{zf'(z)}{f(z)}\Big|<\frac{(1-\alpha)\pi}{2} \quad \mbox{for $ z\in\mathbb{D}$},
$$
where $\alpha\in[0,1].$ He then shows that a domain $\Omega\neq \mathbb{C}$ is
$\alpha$--accessible if and only if $\Omega=f(\mathbb{D}),$ where $f\in {\mathcal S}_{1-\alpha}.$ Note that ${\mathcal S}_{1}$ is the
usual class of all starlike functions.

In \cite{LSt}, Liczberski and Starkov generalized  this result for biholomorphic function $f$ defined on the open unit Euclidean
ball $\mathbb{D}^N,$ $N\geq1,$  such that $f(0)=0.$ They showed that $\alpha$--accessibility is a hereditary property for
such functions $f$, i.e. if $r\in(0,1),$ then $f(r \mathbb{D}^N)$ is also an $\alpha$--accessible domain.

It is known that convexity and starlikeness are not hereditary properties for harmonic functions $f:\,\mathbb{D}\rightarrow\mathbb{C},$ $f(0)=0$
(see \cite{ChDOs}). Therefore, it is natural to consider a subclass of harmonic functions $f$ such
that $f(r\mathbb D)$ is a domain that is  starlike w.r.t the origin (resp. convex) for all $r\in(0,1).$
Such functions were introduced and studied in \cite{ChDOs}.
Now, we will consider the class of $p$--harmonic functions such that $\alpha$--accessibility is a hereditary property in this general setting.

\bdefe \label{def1}
A $p$--harmonic function $\Phi$ in $\mathbb{D}$ with $\Phi(0)=0$
is said to be fully $\alpha$--accessible, $\alpha\in[0,1),$ if for
every $r\in(0,1]$, the function $\Phi$ maps the disk $r\mathbb{D}$ onto an $\alpha$--accessible domain with respect to the origin.
\edefe

\bthm \label{th1}
Let $\Phi(z)=\sum\limits_{k=1}^p|z|^{2(k-1)}F_{p-k+1}(z)$ be a $p$-harmonic function in $\mathbb{D}$
with $\Phi(0)=0$, where $F_{p-k+1}=h_{p-k+1}+\overline{g}_{p-k+1}$ $(k=1,\dots,p)$ are harmonic functions in $\mathbb{D},$
$h_{p-k+1}$ and ${g}_{p-k+1}$ are analytic in $\mathbb{D}.$ Suppose $J_{\Phi}(z)>0$ in $\mathbb{D}.$ Then $\Phi$ is fully
$\alpha$-accessible, $\alpha\in[0,1),$ if and only if the following inequality holds for $z\in\mathbb{D}$:
\begin{equation}\label{eq1}
\sum\limits_{k=1}^p|z|^{2(k-1)}{\rm Re}\,\{
z(h'_{p-k+1}\overline{\Phi}-g'_{p-k+1}{\Phi})\}\geq\sin\frac{\alpha\pi}{2}\cdot|\Phi|\cdot L,
\end{equation}
where $L=\sqrt{A^2+B^2}$ with
\be\label{eq-new5}
\left \{ \begin{array}{l}
A ={\rm Re}\left\{z \sum\limits_{k=1}^p|z|^{2(k-1)}(h'_{p-k+1}-g'_{p-k+1})\right\}\\[4mm]
B={\rm Im}\left\{z \sum\limits_{k=1}^p|z|^{2(k-1)}(h'_{p-k+1}+g'_{p-k+1})\right\}.
\end{array}
\right.
\ee
\ethm\bpf
Consider the function
$$ \Psi(w)=\|\Phi^{-1}(w)\|^2-r^2=\Phi^{-1}(w)\cdot\overline{\Phi^{-1}(w)}-r^2, \quad \mbox{$w\in\mathbb{D}$,}
$$
where $z=\Phi^{-1}(w)$ denotes the inverse mapping so that $\Phi^{-1}(\Phi (z))=z.$
Note that $ \Psi(w)<0$ for $w\in\Phi(r\mathbb{D})$ and $\Psi(w)=0$ for $w\in\partial\Phi(r\mathbb{D}).$

By Theorem \Ref{ThB}, $\Phi$ is fully $\alpha$--accessible if and only if
\begin{equation}\label{eq2}
 \left  \langle \frac{\grad \Psi(w)}{\|\grad
\Psi(w)\|}, \displaystyle\frac{w}{\|w\|}\right  \rangle\geq
\sin\displaystyle\frac{\alpha\pi}{2}
\end{equation}
for every $w=u+iv\in\partial\Phi(r\mathbb{D})$ and for all $r\in(0,1).$

Let us compute $\grad \Psi=(\Psi_u,\Psi_v).$ For this, we first calculate
\begin{equation}\label{eq3}
\Psi_u=\left((\Phi^{-1})_w+(\Phi^{-1})_{\overline
w}\right)\cdot\overline{\Phi^{-1}}+\Phi^{-1}\cdot
\left(\overline{(\Phi^{-1})_{\overline
w}}+\overline{(\Phi^{-1})_{w}}\right),
\end{equation}
and
\begin{equation}\label{eq4}
\Psi_v=\left((\Phi^{-1})_w\thinspace i- (\Phi^{-1})_{\overline w}~
i\right)\cdot\overline{\Phi^{-1}}
+\Phi^{-1}\cdot\left(\overline{(\Phi^{-1})_{\overline w}}~
i-\overline{(\Phi^{-1})_{w}}\thinspace i\right).
\end{equation}
Applying the chain rule in the composition mapping $\Phi^{-1}(\Phi (z))=z$ for the derivatives $\partial /\partial z$ and $\partial /\partial\overline{z}$, we obtain the
following formulas (see for instance \cite[p.\thinspace 9]{Al} and Duren \cite{Du}):
$$(\Phi^{-1})_w=\frac{\overline{(\Phi_{z})}}{J_{\Phi}} ~\mbox{ and }~(\Phi^{-1})_{\overline w}=-\frac{{\Phi_{\overline z}}}{J_{\Phi}}.
$$
Using these two formulas, we may rewrite \eqref{eq3} and \eqref{eq4} as
\be\label{eq-new6}
\Psi_u=\left(\frac{\overline{(\Phi_z)}-\Phi_{\overline
z}}{J_{\Phi}}\right)\overline{z}+z\left(\frac{-\overline{({\Phi}_{\overline z})}+\Phi_{z}}{J_{\Phi}}\right)
=\frac{2}{J_{\Phi}}{\rm Re}\left\{z\left (\Phi_{z}-\overline{(\Phi_{\overline{z}})}\right) \right\},
\ee
and
\be\label{eq-new6a}
 \Psi_v=i\left(\frac{\overline{(\Phi_z)}+\Phi_{\overline
z}}{J_{\Phi}}\right)\overline{z}+zi\left(\frac{-\overline{(\Phi_z)}-\Phi_{\overline{z}}}{J_{\Phi}}\right)
=\frac{2}{J_{\Phi}}{\rm Im}\left\{z\left(\Phi_{z}+\overline{(\Phi_{\overline{z}})}\right)\right\}.
\ee
Using $ |z|^{2(k-1)} = z^{k-1}({\overline z})^{k-1}$, we consider the given $p$-harmonic function
in the form
$$\Phi(z)=F_{p}(z)+\sum\limits_{k=2}^p z^{k-1}({\overline z})^{k-1} F_{p-k+1}(z).
$$
 Then we can easily see that
$$\Phi_{z}
=h'_{p}+\overline z\sum\limits_{k=2}^p
(k-1)|z|^{2(k-2)}F_{p-k+1}+\sum\limits_{k=2}^p
|z|^{2(k-1)}h'_{p-k+1},
$$
and similarly,
$$\overline{\Phi _{\overline{z}}}
=g'_{p}+\overline z\sum\limits_{k=2}^p(k-1)|z|^{2(k-2)}\overline
F_{p-k+1}+\sum\limits_{k=2}^p |z|^{2(k-1)}g'_{p-k+1}.
$$
Using the last two relations and \eqref{eq-new6}, we find that
$$\Psi_{u}
=\frac{2}{J_{\Phi}}{\rm Re}\left\{z\sum\limits_{k=1}^p |z|^{2(k-1)}(h'_{p-k+1}-g'_{p-k+1})\right\} =\frac{2}{J_{\Phi}}A.
$$
Similarly, using these and \eqref{eq-new6a}, we have
$$\Psi_{v}
=\frac{2}{J_{\Phi}}{\rm Im}\left\{z\sum\limits_{k=1}^p
|z|^{2(k-1)}(h'_{p-k+1}+g'_{p-k+1})\right\} =\frac{2}{J_{\Phi}}B.
$$
Here $A$ and $B$ are given by \eqref{eq-new5}. Now, for $w=({\rm Re}\,\Phi,{\rm Im}\,\Phi)$ we calculate the inner product
\beqq
(\grad \Psi, w)& =&\Psi_u\cdot u+\Psi_v\cdot v\\
&=&\frac{2}{J_{\Phi}}\sum\limits_{k=1}^p |z|^{2(k-1)}\Big ({\rm Re}\left\{z(h'_{p-k+1}-g'_{p-k+1})\right\}{\rm Re}\,\Phi \\
&& \hspace{2cm}+{\rm Im}\left\{z(h'_{p-k+1}+g'_{p-k+1})\right\}{\rm Im}\,\Phi\Big )\\
&=&\frac{2}{J_{\Phi}}\sum\limits_{k=1}^p |z|^{2(k-1)}{\rm Re}\left\{z h'_{p-k+1}\overline{\Phi}-z g'_{p-k+1}\Phi\right\}
\eeqq
and thus, \eqref{eq2} is seen to be equivalent to \eqref{eq1}. The proof is complete.
\epf

In the case $\alpha=0$ of Theorem \ref{th1} we get the criterion of fully starlikeness for $p$--harmonic functions.

\bcor
Let $\Phi(z)$ be a $p$--harmonic function defined as in Theorem \ref{th1}. Then $\Phi$ is
fully starlike if and only if
$$\sum\limits_{k=1}^p |z|^{2(k-1)}{\rm Re}\left\{z h'_{p-k+1}(z)\overline{\Phi(z)}\right\}\geq\sum\limits_{k=1}^p
|z|^{2(k-1)}{\rm Re}\left\{z g'_{p-k+1}(z)\Phi(z)\right\}, \quad  z\in\mathbb{D}.
$$
\ecor

The case $p=2$ of Theorem \ref{th1} gives the following.

\bcor \label{cor0}
Let $\Phi(z)=|z|^2F_1(z)+F_2(z),$ $\Phi(0)=0,$ be a sense-preserving biharmonic function in $\mathbb{D},$ where
$F_k=h_k+\overline{g}_k$  $(k=1,2)$ are harmonic in $\mathbb{D}.$ Then $\Phi$ is fully
$\alpha$-accessible with $\alpha\in[0,1)$ if and only if for every $ z\in\mathbb{D}$
\begin{equation}\label{eq*}
 |z|^{2}{\rm Re}\left\{
z(h'_1\overline{\Phi}-g'_1\Phi)\right\}+{\rm Re}\left\{
z(h'_2\overline{\Phi}-g'_2\Phi)\right\}\geq\sin\frac{\alpha\pi}{2}\cdot|\Phi|\cdot
L,
\end{equation}
where
\beqq
L^2&=&|z|^{6}(|h'_1|^{2}+|g'_1|^{2})+|z|^{2}(|h'_2|^{2}+|g'_2|^{2})+2|z|^{4}{\rm Re}\left\{h'_1{\overline{h}_2}'+g'_1{\overline{g}_2}'\right\}-\\
&& \hspace{1cm}
-2{\rm Re}\left\{z^{2}(h'_2g'_2+|z|^{4}h'_1g'_1)\right\}-2|z|^{2}{\rm Re}\left\{z^{2}(h'_1g'_2+g'_1h'_2)\right\}.
\eeqq
\ecor
\bpf
Set $L=\sqrt{A^2+B^2}$, where $A$ and $B$ are defined as in \eqref{eq-new5} with $p=2$. Then, $L$ may be computed explicitly so that
\beqq
L^2&=&\big ({\rm Re\,}\left\{z\thinspace|z|^{2}(h'_1-g'_1)+z(h'_2-g'_2)\right\}\big )^2
+\big ({\rm Im\,}\left\{z\thinspace|z|^{2}(h'_1+g'_1)+z(h'_2+g'_2)\right\}\big )^2\\
&=&\left(|z|^{2}{\rm Re}\,\{z h'_1\}-|z|^{2}{\rm Re}\,\{zg'_1\}+{\rm Re}\,\{z h'_2\}-{\rm Re}\,\{z g'_2\}\right)^2\\
&& +\left(|z|^{2}{\rm Im}\,\{z h'_1\}+|z|^{2}{\rm Im}\,\{zg'_1\}+{\rm Im}\,\{z h'_2\}+{\rm Im}\,\{z g'_2\}\right)^2\\
&=&|z|^{6}(|h'_1|^{2}+|g'_1|^{2})+|z|^{2}(|h'_2|^{2}+|g'_2|^{2})\\
&& +2|z|^{4} \left(-{\rm Re}\,\{z h'_1\}{\rm Re}\,\{zg'_1\}+{\rm Im}\,\{z h'_1\}{\rm Im}\,\{zg'_1\}\right)\\
&&+2|z|^{2} \left({\rm Re}\,\{z h'_1\}{\rm Re}\,\{zh'_2\}+{\rm Im}\,\{z h'_1\}{\rm Im}\,\{zh'_2\}\right)\\
&&+ 2|z|^{2}\left(-{\rm Re}\,\{z h'_1\}{\rm Re}\,\{zg'_2\}+{\rm Im}\,\{z h'_1\}{\rm Im}\,\{zg'_2\}\right)\\
&&+2|z|^{2}\left(-{\rm Re}\,\{z g'_1\}{\rm Re}\,\{zh'_2\}+{\rm Im}\,\{z h'_1\}{\rm Im}\,\{zh'_2\}\right)\\
&&+2|z|^{2}\left({\rm Re}\,\{z g'_1\}{\rm Re}\,\{zg'_2\}+{\rm Im}\,\{z g'_1\}{\rm Im}\,\{zg'_2\}\right)\\
&&+2\left(-{\rm Re}\,\{z h'_2\}{\rm Re}\,\{zg'_2\}+{\rm Im}\,\{z h'_2\}{\rm Im}\,\{zg'_2\}\right)\\
&=&|z|^{6}(|h'_1|^{2}+|g'_1|^{2})+|z|^{2}(|h'_2|^{2}+|g'_2|^{2})-2|z|^{4} {\rm Re}\,\{z^2 h'_1 g'_1\}+ 2|z|^{4}{\rm Re}\,\{h'_1 {\overline{h}_2}'\}\\
&&-2|z|^{2} {\rm Re}\,\{z^2 h'_1 g'_2\}-2|z|^{2} {\rm Re}\,\{z^2 g'_1 h'_2\}+2|z|^{4}{\rm Re}\,\{g'_1 {\overline{g}_2}'\}-2{\rm Re}\,\{z^2 h'_2 g'_2\}.
\eeqq
Simplifying the last expression gives the inequality \eqref{eq*}.
\epf

\bcor \label{cor2}
Let $\Phi(z)=h(z)+\overline{g(z)}$ be a harmonic function in $\mathbb{D},$ $\Phi(0)=0,$ $J_{\Phi}(z)>0$ in
$\mathbb{D}.$ Then $\Phi$ is fully $\alpha$-accessible with $\alpha\in[0,1)$ if and only if
$$|h(z)|^{2}\cdot{\rm Re}\left\{\frac{zh'(z)}{h(z)}\right\}-|g(z)|^{2}\cdot{\rm Re}\left\{\frac{zg'(z)}{g(z)}\right\}-{\rm Re}\left\{
z(g'(z)h(z)-h'(z)g(z))\right\}$$
$$\geq\sin\frac{\alpha\pi}{2} \left |zh'(z)-\overline{zg'(z)}\right | \, \left |h(z) +\overline{g(z)}\right |
$$
for $z\in\mathbb{D}$.
\ecor \bpf
Suppose that $\Phi=h+\overline{g}$ is harmonic. Then we may set $p=1$ in Theorem \ref{th1} and thus,
the left-hand side of inequality \eqref{eq1} equals
$${\rm Re}\left\{ z(h'(\overline h+g)-g'(h+\overline g))\right\}
=|h|^{2}{\rm Re}\left\{\frac{zh'}{h}\right\}-|g|^{2}{\rm Re}\left\{\frac{zg'}{g}\right\}-{\rm Re}\left\{
z(g'h-h'g)\right\}.
$$
Moreover, $L^2$  can be written in the form
\beqq
L^2&=& \big ({\rm Re}\left\{z (h'-g')\right\}\big )^2+ \big ({\rm Im}\left\{z (h'+g')\right\}\big )^2\\
&=&|z|^{2}(|h'|^{2}+|g'|^{2})-2{\rm Re}\left\{z^2 h'g'\right\}\\
&=&|zh'-\overline{zg'}|^{2}.
\eeqq
The proof is complete.
\epf

The case $\alpha=0$ of Corollary \ref{cor2} reduces to the well-known result of Chuaqui et al. \cite{ChDOs} for
fully starlikeness of harmonic functions.
Moreover, if the function $\Phi(z)$ in Corollary \ref{cor2} is analytic, then we have the following well-known result.

\bcor
An analytic function $\Phi(z)$ in $\ID$ with $\Phi (0)=0$ is fully $\alpha$-accessible with $\alpha\in[0,1)$ if and only if
$${\rm Re}\left\{\frac{z\thinspace\Phi'(z)}{\Phi(z)}\right\}\geq\sin\frac{\alpha\pi}{2}\left|\frac{z\thinspace\Phi'(z)}{\Phi(z)}\right|
~\mbox{ for $z\in\mathbb{D}$}.
$$
\ecor

The last inequality is equivalent to the inequality
\begin{equation}\label{eq5}
\left|\arg\left(\frac{z\thinspace\Phi'(z)}{\Phi(z)}\right)\right|\leq\frac{\pi}{2}(1-\alpha).
\end{equation}
The condition \eqref{eq5} is known, because $\Phi(\mathbb{D})$ is $\alpha$-accessible ($\alpha\in[0,1))$
if and only if $\Phi$ is fully $\alpha$-accessible (cf. \cite{LSt}), and $\Phi(\mathbb{D})$
is $\alpha$-accessible which in turn is equivalent to the condition \eqref{eq5} (cf. \cite{St1, St2, BrK, Sug}).

Clearly, the case $\alpha =0$ of \eqref{eq5} is the condition for the
starlikeness of the analytic function $\Phi(z)$ in $\ID$.


Next, we let $k=p$, $F_{p-k+1}(z)\equiv 0$ for $k=1,\ldots, p-1$ and $F_1(z)$ be analytic in $\ID$. Then the corresponding $\Phi(z)$ in
Theorem \ref{th1} reduces to
$$\Phi(z)=|z|^{2(p-1)}F_1(z)
$$
and Jacobian $J_\Phi$ takes an easy form and hence (see also \cite[Corollary 1.2]{Li})
\bee
\item[(i)] $J_\Phi (z)=0$ if and only if $\ds \left |  \frac{zF_1'(z)}{F_1(z)} + p-1 \right |  = p-1$ or $z=0$;
\item[(ii)] $J_\Phi(z) >0$ if and only if $\ds \left |  \frac{zF_1'(z)}{F_1(z)} + p-1 \right |  > p-1$ and  $z\neq 0$;
\item[(iii)] $J_\Phi (z)<0$ if and only if $\ds \left |  \frac{zF_1'(z)}{F_1(z)} + p-1 \right |  < p-1$ and  $z\neq 0$.
\eee
It was shown in \cite[Corollary 1.2]{Li} that $\Phi$ of the above form is starlike and univalent in $\ID$ whenever $F_1$ is starlike in $\ID$.
On the other hand, using Theorem \ref{th1} this result can be generalized as follows.

\bcor\label{cor-new1} Let $\Phi(z)=|z|^{2(p-1)}F_1(z)$, where $F_1$
is analytic and  $\alpha$-accessible in $\ID$. Then $\Phi$ is fully
$\alpha$-accessible and univalent in $\ID$. \ecor \bpf
A computation shows that the inequality \eqref{eq1} reduces to
$${\rm Re}\left\{\frac{zF_1'(z)}{F_1(z)}\right\}\geq \left|\frac{z\thinspace F_1'(z)}{F_1(z)}\right| \sin\frac{\alpha\pi}{2}
~\mbox{ for $z\in\mathbb{D}$},
$$
and the desired conclusion follows from Theorem \ref{th1}.
\epf

For simplicity, we now present a couple of examples associated with Corollary \ref{cor-new1}.

\begin{example}
Let $F_1(z)=z+\lambda z^n$, where $0<|\lambda |\leq 1/n$ and $n\geq 2$. Then we have
$$\left |\frac{zF_1'(z)}{F_1(z)} -1\right |
=  \left |\frac{\lambda (n-1)z^{n-1}}{1+\lambda z^{n-1}} \right |
< \frac{|\lambda | (n-1)}{1-|\lambda | }
$$
which shows that
$$\left|\arg\left(\frac{z F_1'(z)}{F_1(z)}\right)\right|\leq\arcsin \left (
\frac{|\lambda | (n-1)}{1-|\lambda | } \right ).
$$
Consequently, for $0<|\lambda |\leq 1/n$, the function
$\Phi(z)=|z|^{2(p-1)}(z+\lambda z^n)$ is fully $\alpha$-accessible
and univalent in $\ID$ with
$$\alpha:= \alpha (n,\lambda) = 1-\frac{2}{\pi}\arcsin \left ( \frac{|\lambda | (n-1)}{1-|\lambda | } \right ).
$$
The graph of $\Phi$ with $\lambda =1/n$ and various values of $p$ and $n$ are shown in Figure \ref{fig-F332}. Observe that
$\alpha (n,1/n)=1/2$ and thus,
$\Phi(z)=|z|^{2(p-1)}(z+(1/n) z^n)$ is fully $\alpha$-accessible and univalent in $\ID$ with $\alpha =1/2$.

\begin{figure}[H]
\begin{center}
\includegraphics[height=5.5cm, width=5.5cm, scale=1.3]{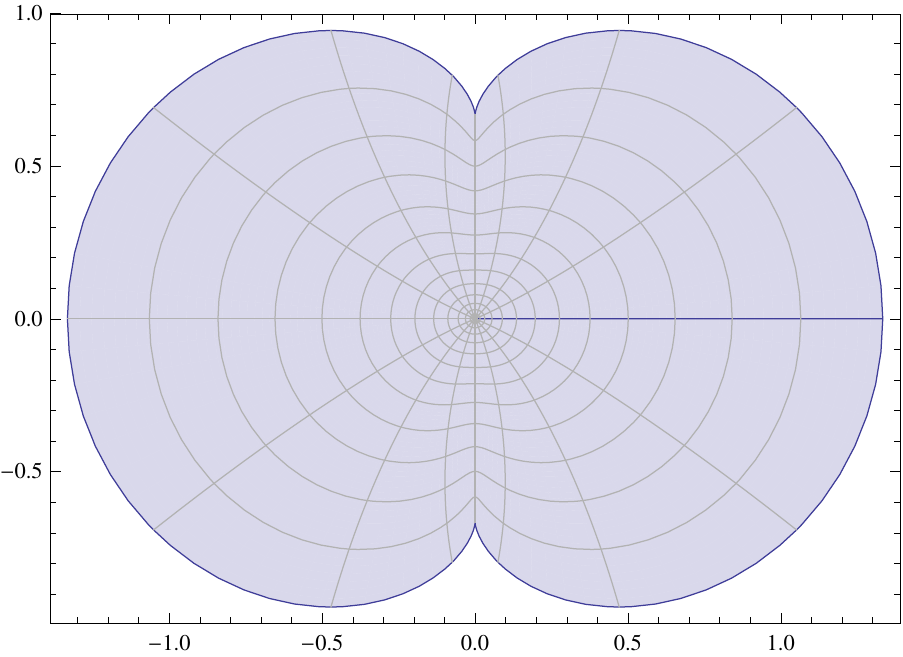}
\hspace{1cm}
\includegraphics[height=5.5cm, width=5.5cm, scale=1.3]{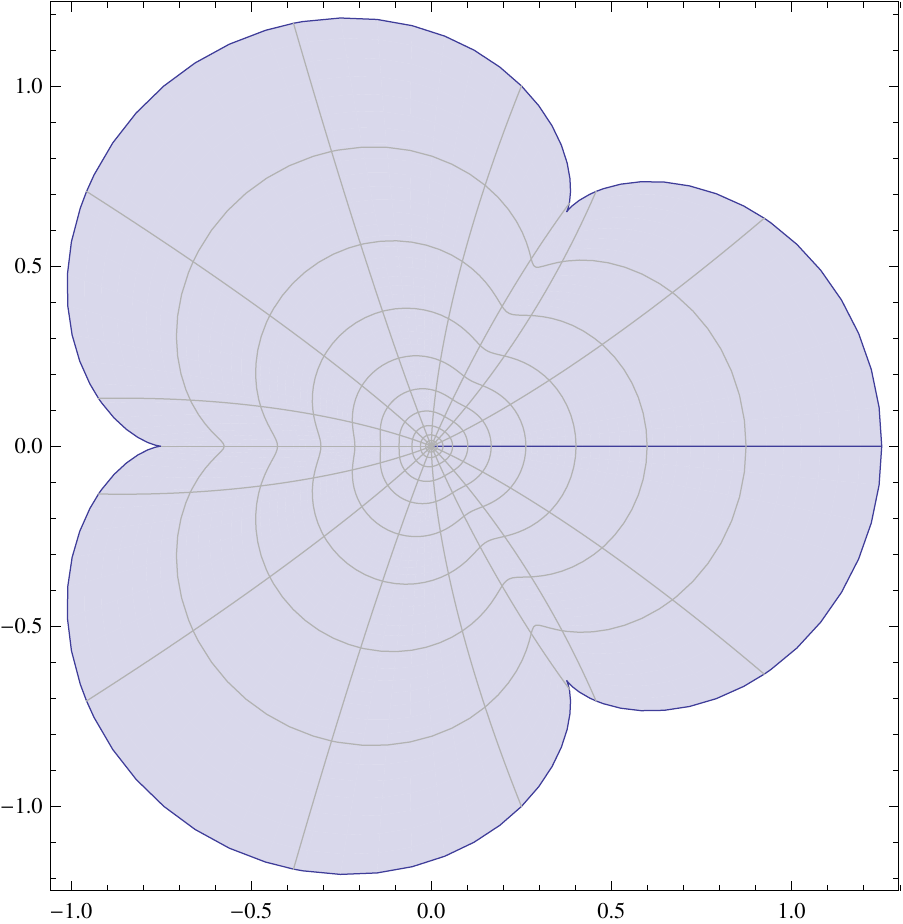}
\end{center}
(1) $n=3$, $p=2$ \hspace{4.5cm} (2) $n=4$, $p=3$

\vspace{0.25cm}

\begin{center}
\includegraphics[height=5.5cm, width=5.5cm, scale=1]{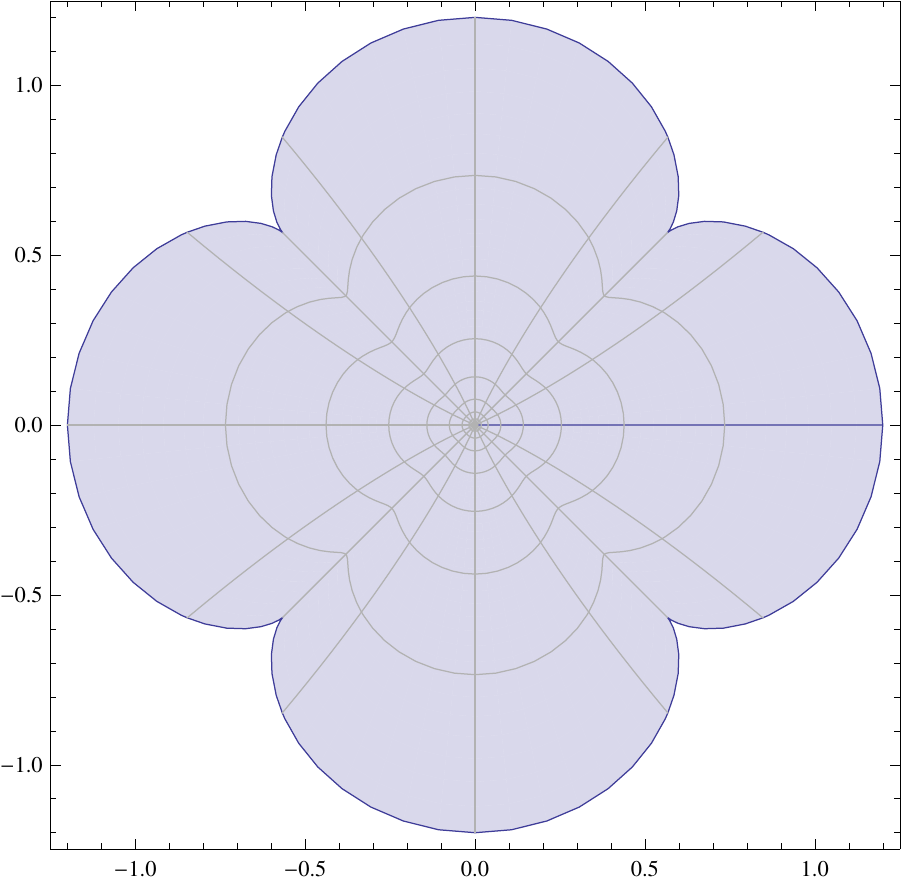}
\hspace{1cm}
\includegraphics[height=5.5cm, width=5.5cm, scale=1]{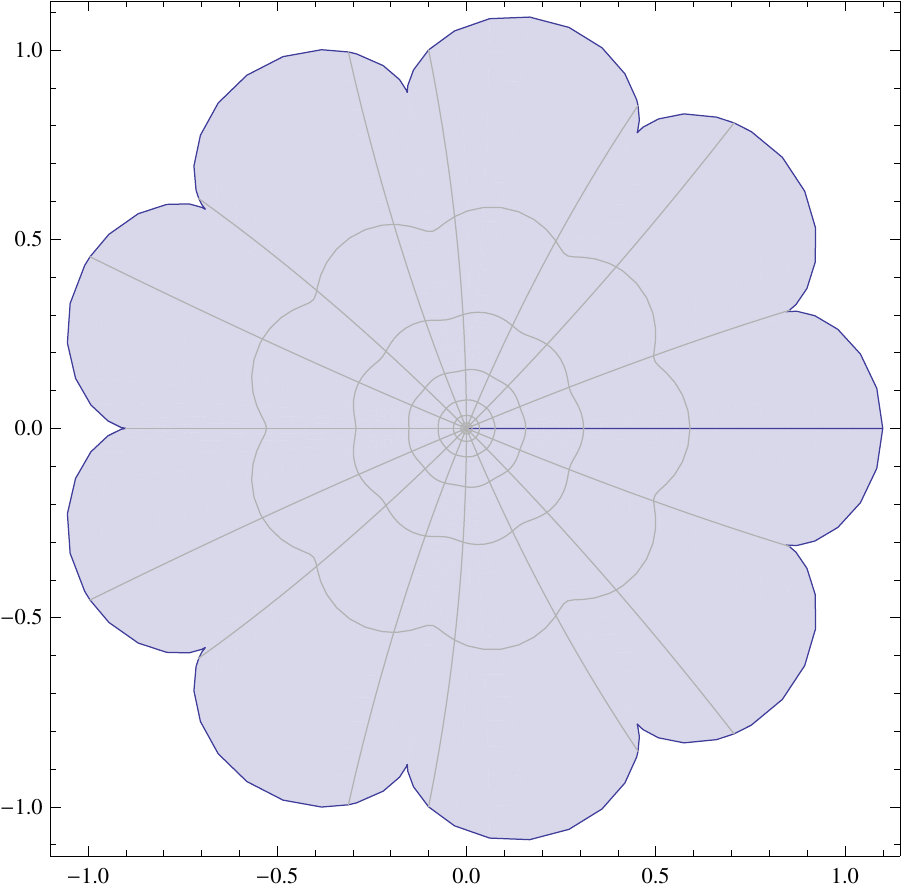}
\end{center}
(3) $n=5$, $p=4$ \hspace{4.5cm} (4) $n=10$, $p=5$

\caption{The range of $\Phi(z)=|z|^{2(p-1)}(z+(1/n)z^n)$ for certain values of $n$ and $p$.\label{fig-F332}}
\end{figure}
\end{example}

\begin{example}
Consider the $p$-harmonic function $F(z)$ in $\mathbb D$ defined by
$$F(z)=z-\lambda |z|^{2(p-1)}:=|z|^{2(p-1)}F_1(z)+F_p(z),
$$
where $\lambda \in\IC$ with $0<|\lambda|<\frac{1}{2(p-1)}$ and $p\geq 2$. It follows that
\beqq
J_F(z)&=&\left|1-\lambda(p-1)\thinspace|z|^{2(p-2)}\thinspace\overline z\right|^2-\left|\lambda(p-1)\thinspace|z|^{2(p-2)} z\right|^2\\
&=& 1-2 (p-1)\thinspace|z|^{2(p-2)}\thinspace{\rm Re}\,\{\overline{\lambda} z\}\\
&\geq &  1-2|\lambda| (p-1)>0
\eeqq
which shows that $F$ is sense-preserving in $\ID$. We now apply Theorem \ref{th1} to the function $F(z).$ The corresponding
inequality
$${\rm Re}\left\{z\thinspace(\overline z-\overline{\lambda}\thinspace
|z|^{2(p-1)})\right\}\geq\sin\frac{\alpha\pi}{2}\thinspace\left|z-\lambda\thinspace
|z|^{2(p-1)}\right|\cdot\left( ({\rm Re}\, z)^2+({\rm Im}\, z)^2\right)^{1/2}
$$
may be rewritten as
$$ {\rm Re}\left\{1-\overline{\lambda} \thinspace |z|^{2(p-2)}\thinspace
z\right\}\geq\sin\frac{\alpha\pi}{2}\thinspace\left|1-\lambda\thinspace
|z|^{2(p-2)}\overline{z}\right|.
$$
From this it follows that $F(z)$ is fully $\alpha$--accessible, where
$$\sin\frac{\alpha\pi}{2}=\min\limits_{z\in{\overline{\mathbb
D}}}\left (\frac{1-\thinspace |z|^{2(p-2)}\thinspace{\rm Re}\,\{
\overline{\lambda} z\}}{\left|1-\thinspace
\lambda |z|^{2(p-2)}\overline{ z}\right|}\right )\geq\frac{1-|\lambda|}{1+|\lambda|}.
$$
Thus, we can choose $\alpha$ with
$$ \alpha\geq\frac{2}{\pi}\arcsin\frac{1-|\lambda|}{1+|\lambda|}=:\alpha_0
$$
which shows that $F(z)$ is fully $\alpha_0$--accessible.
\end{example}

\medskip

\subsection*{Acknowledgements} The authors thank  Sh. Chen, S. Yu. Graf, Sairam Kaliraj, and  V. V. Starkov for fruitful discussion
and comments on this manuscript. The research was supported by the project RUS/RFBR/P-163 under Department of Science \& Technology (India).
The first and the second authors were supported by Strategic Development Program of Petrozavodsk State University. The second author was supported by a
grant from the Simons Foundation and RFBR (project N 14-01-00510a). The third author is currently on leave from Indian Institute of
Technology Madras.


\end{document}